\documentclass[12pt]{article}
\usepackage{amsmath,amsthm, amssymb, amsfonts,latexsym,amscd}
\usepackage[dvips]{graphicx}
\newtheorem{theorem}{Theorem}[section]
\newtheorem{corollary}[theorem]{Corollary}
\newtheorem{lemma}[theorem]{Lemma}
\newtheorem{proposition}[theorem]{Proposition}

\newtheorem{example}[theorem]{Example}
\begin{document}

\title{Minimizing Laplacian spectral radius of unicyclic graphs with fixed girth}
\author{Kamal Lochan Patra\footnote{klpatra@niser.ac.in}\; and Binod Kumar Sahoo\footnote{bksahoo@niser.ac.in}}
\date{}
\maketitle
\begin{center}
School of Mathematical Sciences\\
National Institute of Science Education and Research\\
Sainik School Post, Bhubaneswar--751005, India.\\
\end{center}

\begin{abstract}
In this paper we consider the following problem: Over the class of all simple connected unicyclic graphs on $n$ vertices with girth $g$ ($n,g$ being fixed), which graph minimizes the Laplacian spectral radius? We prove that the graph $U_{n,g}$ (defined in Section 1) uniquely minimizes the Laplacian spectral radius for $n\geq 2g-1$ when $g$ is even and for $n\geq 3g-1$ when $g$ is odd.\\

\noindent Keywords: Laplacian matrix; Laplacian spectral radius; girth; unicyclic graph.
\end{abstract}

\section{Introduction}

Let $G=(V, E)$ be a simple graph with vertex set $V=\{v_1, v_2,
\ldots, v_n\}$ and edge set $E.$ The {\em adjacency matrix} of $G$ is defined as $A(G)=(a_{ij}),$ where $a_{ij}$ is
equal to $1$ if $\{v_i,v_j\} \in E$ and $0$ otherwise. Let $D(G)$ be
the {\em diagonal matrix} of $G$ whose $i$-th diagonal entry is the
degree of the vertex $v_i$ of $G.$ The {\em Laplacian matrix} $L(G)$ of $G$ is defined by: $L(G)=D(G)-A(G)$. Clearly $L(G)$ is real symmetric. It is known that $L(G)$
is a positive semi-definite matrix. So all its eigenvalues are real and non-negative. Since the sum of the
entries in each row of $L(G)$ is zero, the all
one vector ${\bf e}=[1,\cdots,1]^{T}$ is an eigenvector of $L(G)$ corresponding to the smallest eigenvalue zero.
Here $X^T$ denotes the transpose of a given matrix $X$. For more about the Laplacian matrix and its eigenvalues we refer the reader to \cite{merris,merris1,mohar} and the
references therein.

The largest eigenvalue of $L(G)$ is called the {\em Laplacian spectral radius} of $G$, we denote it by $\lambda(G)$. Among all trees on $n$ vertices, the
Laplacian spectral radius is uniquely minimized by the path, and uniquely maximized by the star.  The tree
that uniquely maximizes the Laplacian spectral radius over all trees on $n$ vertices with fixed diameter (respectively, with fixed number of pendant vertices) is characterized in \cite{lal} (respectively, in \cite{guo}). Over all unicyclic graphs on $n$ vertices, the cycle has the minimum Laplacian spectral radius and the graph obtained by joining $n-3$ pendant vertices to a vertex of a 3-cycle has the maximum Laplacian spctral radius.

The second smallest eigenvalue of $L(G)$ is called the {\it algebraic connectivity} of $G$ \cite{fiedler}. The Laplacian spectral radius of $G$ is related to the algebraic connectivity of the complement graph of $G$. The study of the algebraic connectivity and the Laplacian spectral radius of graphs has received a good deal of attention in the recent past. As far as the class of connected unicyclic graphs on $n$ vertices with fixed girth is concerned, the problem of minimizing and maximizing the algebraic connectivity has been studied in \cite{FKP1} and \cite{FKP2}, respectively, and that of maximizing the Laplacian spectral radius is done by Guo (\cite{guo1}, Corollary 4.1). However, the question of minimizing the Laplacian spectral radius has not been studied so far and we deal with this problem in the present paper.

Recall that the {\em girth} of a graph $G$ is the length of a shortest cycle (if any) of
$G$. A graph $G$ is {\em unicyclic} if it has exactly one cycle, that is, the number of edges and the number of vertices of $G$ are the same.
Consider a cycle on $g\geq 3$ vertices and append a pendent vertex of the path on $n-g\;(n> g)$
vertices to a vertex of the cycle. The new graph thus obtained is a unicyclic graph on $n$ vertices
with girth $g$. We denote it by $U_{n,g}$ (see Figure \ref{fig:uni}). In this paper, we prove

\begin{theorem}\label{main} Let $G$ be a unicyclic graph on $n$ vertices with girth $g$ which is not isomorphic to $U_{n,g}$. Then the following hold.
\begin{enumerate}
\item[(1)] If $g$ is even and $n\geq 2g-1$, then $\lambda(U_{n,g})< \lambda(G)$.

\item[(2)] If $g$ is odd and $n\geq 3g-1$, then $\lambda(U_{n,g})< \lambda(G)$.
\end{enumerate}
\end{theorem}

In Section 2, we recall some basic definitions and results from the literature which are needed in the subsequent sections. In Section 3, we study the Laplacian spectral radius of the graph $U_{n,g}$. Finally, we prove Theorem \ref{main} in Section 4.

\section{Preliminaries}

All graphs considered in this paper are finite, simple and connected. For any graph $G$, define $B(G)=D(G)+A(G)$. We denote by $\mu(G)$ the largest eigenvalue of $B(G)$. By {\em Rayleigh-Ritz} theorem (\cite{horn}, p.176), we have
$$\lambda(G)=\max_{X\in W}\;X^TL(G)X;$$
$$\mu(G)=\max_{X\in W}\;X^TB(G)X;$$
where $W=\{X\in\mathbb{R}^n\;|\;X^TX=1\}$.
Since $G$ is connected, $B(G)$ is a nonnegative irreducible matrix. So $\mu(G)$ is simple and there exists a positive eigenvector of $B(G)$ corresponding to $\mu(G)$. This is a consequence of the Perron-Frobenius theory. If $X$ is a unit eigenvector of $B(G)$ corresponding to $\mu(G)$, then we have
$$\mu(G)=X^TB(G)X=\sum_{\{v_i, v_j\}\in E}(x_i + x_j)^2,\;\mbox{where}\;\;X^T=[x_1,x_2,\cdots,x_n].$$

Recall that $G$ is {\em bipartite} if its vertex set $V$ is a disjoint union of two sets $V_1$ and $V_2$ such that every edge in $G$ joins a vertex of $V_1$ to a vertex of $V_2$. The next proposition relates $L(G)$ and $B(G)$ for a bipartite graph $G$.

\begin{lemma}[\cite{grone1}, p.220]\label{lem:unitary}
Let $G$ be a bipartite graph. Then $B(G)$ and $L(G)$ are unitarily similar. In
particular, $\lambda(G)$ is simple.
\end{lemma}

By Lemma \ref{lem:unitary} and the Perron-Frobenius theory, the following is immediate.

\begin{lemma}\label{lem:bipartite}
Let $G$ be a bipartite graph and $G'$ be a graph obtained from $G$ by adjoining a new vertex to a vertex of $G$.
Then $\lambda(G')>\lambda(G).$
\end{lemma}

\begin{lemma}[\cite{grone1}, p.233]\label{lem:courant}
Let $G$ be a graph and $G'$ be the graph obtained from $G$ by joining two non-adjacent vertices of $G$ with an edge. Then $\lambda(G')\geq\lambda(G)$.
\end{lemma}

\begin{lemma}[\cite{grone}, p.224]\label{lem:maxdegree}
 Let $G$ be a graph on $n\geq 2$ vertices. Then $\lambda(G)\geq \Delta(G)+1,$ where $\Delta(G)$ is the maximum vertex degree of $G$. Further, equality holds if and only if $\Delta(G)=n-1$.
\end{lemma}

Lemma \ref{lem:courant} says that introducing a new edge in a graph $G$ can not decrease the Laplacian spectral
radius. The next result mentions a case in which the Laplacian spectral radius remains the same even after introducing a new edge in
the given graph.

\begin{proposition}[\cite{guo1}, p.712]\label{thm:guo}
 Let $G$ be a graph on $n\geq 2$ vertices and $v$ be a vertex of $G$. Let $G_s^k$ be the graph obtained from $G$ by attaching $s\geq 2$ new paths $P_i:vv_{ik}v_{i(k-1)}\cdots v_{i1}\;\;(1\leq i\leq s)$ at $v$ each of length $k\geq 1$.
Let $G_{s,t}^k$ be the graph obtained from $G_s^k$ by adding $t\;(1\leq t\leq\frac{s(s-1)}{2})$ edges among the vertices $v_{11},v_{21},\cdots,v_{s1}$.
Then $\lambda(G_s^k)=\lambda(G_{s,t}^k)$.
\end{proposition}

\begin{figure}[h]
\centering
\includegraphics[scale=.9]{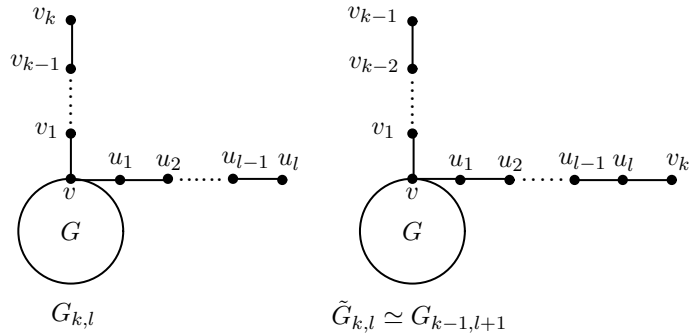}
\caption{Grafting an edge}\label{fig:gkl}
\end{figure}

Let  $G$ be a graph on $n\ge 2$ vertices and $v$ be a vertex of $G$. For $l\geq k \geq 1$, we construct a new graph $G_{k,l}$ from $G$ by attaching two new paths
$P:vv_{1}v_{2}\ldots v_{k}$ and $Q:vu_{1}u_{2}\ldots u_{l}$ of
lengths $k$ and $l$, respectively, at $v$. Let $\widetilde{G}_{k,l}$ be the graph obtained
from $G_{k,l}$  by removing the edge $\{v_{k-1},v_{k}\}$ and adding
the edge $\{u_{l},v_{k}\}$ (see Figure \ref{fig:gkl}). We say that
$\widetilde{G}_{k,l}$ is obtained from $G_{k,l}$ by {\em grafting} an
edge. The next result compares the Laplacian spectral radius of
$G_{k,l}$ and $\widetilde{G}_{k,l}\simeq G_{k-1, l+1}.$

\begin{proposition}[\cite{guo}, p.65,68]\label{thm:gkl}
Let $G$ be a graph on $n\geq 2$ vertices and
$v$ be a vertex of $G$. For $l\geq k\geq 1$, let $G_{k,l}$ be the graph defined as
above.  Then $$\lambda(G_{k-1,l+1})\leq
\lambda(G_{k,l}),$$ with equality if and only if there exists an
eigenvector of $G_{k,l}$ corresponding to $\lambda(G_{k,l})$ whose $v$-th component is zero. In particular, if $G$ is bipartite, then $\lambda(G_{k-1,l+1})< \lambda(G_{k,l})$.
\end{proposition}

\begin{lemma}\label{lem:fx}
Let $f_1(y)=y-1$ and define $f_i(y)=y-2-\frac{1}{f_{i-1}(y)}$ for $i\geq 2$. Then for $i,j\geq 1$,
\begin{enumerate}
\item[$(i)$] $f_i(y)>\frac{y}{y-2}$ for $y\geq 4.383$.
\item[$(ii)$] $f_i(y)> f_{i+1}(y)>1$ for $y\geq 4$.
\item[$(iii)$] $f_i(y)f_{i+1}(y)>f_j(y)$ for $y\geq 4.383$.
\end{enumerate}
\end{lemma}

\begin{proof} The proof of $(i)$ and $(ii)$ is similar to that of Lemma 3.2 in \cite{guo}. We prove $(iii)$ now. We have $f_i(y)f_{i+1}(y)=f_i(y)(y-2-\frac{1}{f_i(y)})=(y-2)f_i(y)-1$. Since $y\geq 4.383$, $(i)$ implies that $f_i(y)f_{i+1}(y) >(y-2)(\frac{y}{y-2})-1=y-1=f_1(y)\geq f_j(y)$. Here the last inequality follows from $(ii)$.
\end{proof}

For an eigenvector $X$ of a graph $G$ corresponding to $\mu(G)$ (or $\lambda(G)$), we associate with $X$ a labeling of $G$ in which a vertex $v_i$ is labeled $x_{v_i}$ or simply $x_i$.

\begin{lemma}\label{lem:fx1}
Let $v$ be a vertex of a bipartite graph $H$ and $G$ be a graph obtained from $H$ by attaching a path $P:v\text{=}v_{0}v_{1}\cdots v_k$ at $v$. Let $X$ be a positive eigenvector of $B(G)$ corresponding to $\mu=\mu(G)$. Then $x_{v_{i}}=f_{k-i}(\mu)x_{v_{i+1}}$ for $0\leq i\leq k-1$, where $f_i(y)$ is the function defined in Lemma \ref{lem:fx}. Further, if $\mu\geq 4$, then $x_{v_i}>x_{v_{i+1}}$ for $0\leq i\leq k-1$.
\end{lemma}

\begin{proof}
From $B(G)X=\mu X$, we have $x_{v_{k-1}}=(\mu-1)x_{v_k}$ and $x_{v_{i-2}}+x_{v_i}=(\mu-2)x_{v_{i-1}}$ for $2\leq i\leq k$. Using these two equations it follows that $x_{v_{i}}=f_{k-i}(\mu)x_{v_{i+1}}$ for $0\leq i\leq k-1$. If $\mu\geq 4$, then $f_{k-i}(\mu)> 1$ by Lemma \ref{lem:fx}$(ii)$. So $x_{v_i}>x_{v_{i+1}}$ for $0\leq i\leq k-1$.
\end{proof}

Let $G$ be a graph with vertex set $V=\{v_1,\cdots,v_n\}$ and Laplacian matrix $L=L(G)=(l_{ij})$. Let $\tau$ be an automorphism of $G$. Since $\tau$ is a permutation of $V$, it induces a permutation matrix $P=(p_{ij})$, where $p_{ij}$ is defined by:
\begin{equation*}
p_{ij}=\left\{\begin{array}{ll}
  1 & \text{ if }\tau(v_j)=v_i \\
  0 & \text{ otherwise }
\end{array}.\right.
\end{equation*}
If $\tau(v_t)=v_i$ and $\tau(v_j)=v_s$, then the $ij$-th entry of $PL$ is $l_{tj}$ and that of $LP$ is $l_{is}$. We have that $v_t\neq v_j$ if and only if $v_i\neq v_s$, and that $v_t$ and $v_j$ are adjacent if and only if $v_i$ and $v_s$ are adjacent. This, together with the fact that $\tau$ preserves the degree of a vertex, implies that $l_{tj}=l_{is}$. Thus $PL=LP$. The same argument implies that $PB=BP$, where $B=B(G)$.

\begin{lemma}\label{lem1:auto}
Let $G$ be a  graph and $\tau$ be an automorphism of $G$. If $\lambda=\lambda(G)$ is a simple eigenvalue of $L=L(G)$ and $X$ is an eigenvector corresponding to $\lambda$, then $|x_{v_k}|=|x_{\tau(v_k)}|$ for every $v_k\in V$.
\end{lemma}

\begin{proof}
We have $LPX=PLX=P\lambda X=\lambda PX$. So $PX$ is also an eigenvector of $L$ corresponding to $\lambda$. Since $\lambda$ has algebraic multiplicity one, $X$ and $PX$ are linearly dependent. So $PX=\alpha X$ for some real number $\alpha$. Since $P$ has finite order (as a group element in $GL(n,\mathbb{R})$), it follows that $\alpha$ is a $k$-th root of unity for some positive integer $k$. So $\alpha=\pm 1$ and $PX=\pm X$. Then $\sum_{j=1}^{n} p_{ij}x_{v_j}=\pm x_{v_i}$. If $\tau(v_k)=v_i$, then $x_{v_k}=\pm x_{v_i}$ and so $|x_{v_k}|=|x_{v_i}|=|x_{\tau(v_k)}|$. This completes the proof.
\end{proof}

\begin{lemma}\label{lem2:auto}
Let $G$ be a bipartite graph and $\tau$ be an automorphism of $G$. Let $X$ be a positive eigenvector of $B(G)$ corresponding to $\mu(G)$. Then $x_{v_k}=x_{\tau(v_k)}$ for $v_k\in V$.
\end{lemma}

\begin{proof}
Since $G$ is bipartite, $\mu(G)$ is a simple eigenvalue of $B(G)$. Now, the proof is similar to that of Lemma \ref{lem1:auto}.
\end{proof}

\section{The Graph $U_{n,g}$}

\begin{figure}[h]
\centering
\includegraphics[scale=.9]{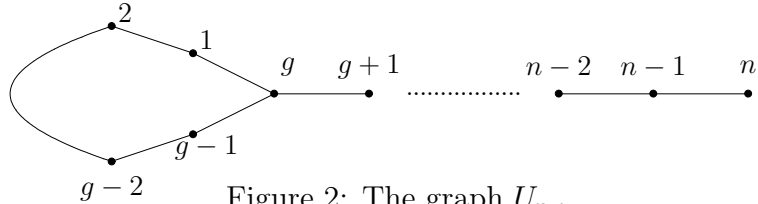}
\caption{The graph $U_{n,g}$}\label{fig:uni}
\end{figure}

In Section 1, we defined the graph $U_{n,g}$ for $n>g$. We take the vertex set $V$ of $U_{n,g}$ as $V=\{1,2,\cdots,n\}$ and the edges of $U_{n,g}$ as shown in Figure \ref{fig:uni}. Note that, by Lemma \ref{lem:maxdegree}, $\lambda(U_{n,g})\geq 4$ with equality if and only if $n=4$.

We first consider the graph $U_{n,g}$ with $g$ even. So $U_{n,g}$ is a bipartite graph and hence $\lambda(U_{n,g})=\mu(U_{n,g})$ by Lemma \ref{lem:unitary}.

\begin{proposition}\label{prop:ung1}
Let $G=U_{n,g}$ with $g$ even and $X=[x_1,\cdots,x_n]^T$ be a positive eigenvector of $B=B(G)$ corresponding to $\mu=\mu(G)$. Then the following hold:
\begin{enumerate}
\item[$(i)$] $x_j=x_{g-j}$ for $j\in\{1,\cdots, \frac{g}{2}\}$.
\item[$(ii)$] $x_j>x_{j+1}$ for $j\in\{g,\cdots, n-1$\}.
\item[$(iii)$] $x_g>x_{1}$ and $x_j>x_{j+1}$ for $j\in\{1,\cdots,\frac{g}{2}-1$\}.
\item[$(iv)$] If $\mu\geq 4.5$, then $x_j>2x_{j+1}$ for $j\in\{g,\cdots,n-1\}$, $x_g<2x_{1}$ and $x_j<2x_{j+1}$ for $j\in\{1,\cdots,\frac{g}{2}-1\}$.
\item[$(v)$] If $\mu\geq 4.5$, then $x_i>x_{g+i}$, where $1\leq i\leq\text{min }\{\frac{g}{2}, n-g\}$.
\end{enumerate}
\end{proposition}

\begin{proof}
$(i)$ The map $\tau$ fixing the vertices $g,g+1,\cdots,n$ and taking $j$ to $g-j$ for $j\in\{1,\cdots,g-1\}$ is an automorphism of $G$. Since $G$ is a bipartite graph,  the result follows from Lemma \ref{lem2:auto}.

$(ii)$ This follows from Lemma \ref{lem:fx1}, since $\mu=\lambda(G)\geq 4$.

$(iii)$ From $BX=\mu X$ at the vertex $\frac{g}{2}$, we have $x_{\frac{g}{2}-1}+2x_{\frac{g}{2}}+x_{\frac{g}{2}+1}=\mu x_{\frac{g}{2}}$. Since $x_{\frac{g}{2}-1}=x_{\frac{g}{2}+1}$ by $(i)$ and $\mu>4$, $x_{\frac{g}{2}}=(\frac{2}{\mu-2})x_{\frac{g}{2}-1}<x_{\frac{g}{2}-1}$. Assume that $x_{j+1}<x_j$ for $j\in\{2,\cdots,\frac{g}{2}-1\}$. Now $x_{j-1}=(\mu-2)x_j-x_{j+1}>(\mu-3)x_j>x_j$. So $x_j>x_{j+1}$ for $j\in\{1,\cdots,\frac{g}{2}$\}. A similar proof holds for $x_g>x_1$.

$(iv)$ For $j=n-1$, we have $x_{n-1}=(\mu-1)x_n>2x_n$. Assume that $x_{j}>2x_{j+1}$ for $j\in\{g+1,\cdots, n-1\}$. Now $x_{j-1}=(\mu-2)x_j-x_{j+1}\geq 2.5x_j-x_{j+1}>2.5x_j-0.5x_j=2x_j$. So $x_j>2x_{j+1}$ for $j\in\{g,\cdots, n-1\}$. From $BX=\mu X$ at the vertex $g$ and using $x_1=x_{g-1}$ by $(i)$, we have $2x_1=(\mu-3)x_g-x_{g+1}\geq 1.5 x_g-x_{g+1}>1.5 x_g-0.5x_g=x_g$. So $x_g<2x_1$. Now at the vertex $1$, $x_{2}=(\mu-2)x_1-x_{g}\geq 2.5 x_1-x_{g}>2.5 x_1-2x_1=\frac{1}{2}x_1$. So $x_1<2x_2$. Assume that $x_{j-1}<2x_{j}$ for $j\in\{2,\cdots,\frac{g}{2}-1\}$. Now $x_{j+1}=(\mu-2)x_j-x_{j-1}\geq 2.5 x_j-x_{j-1}>2.5 x_j-2x_j=\frac{1}{2}x_j$. So $x_j<2x_{j+1}$ for $j\in\{1,\cdots,\frac{g}{2}-1\}$.

$(v)$ By $(iv)$, $2x_1>x_g>2x_{g+1}$. So $x_1>x_{g+1}$. Inductively, for $1\leq i\leq\text{min }\{\frac{g}{2}, n-g\}$, $2x_i>x_{i-1}>x_{g+i-1}>2x_{g+i}$ gives that $x_i>x_{g+i}$.
\end{proof}

\begin{lemma}\label{lem:4.5}
Let $g$ be even, $k=\frac{g}{2}$ and $n=g+k$. Then $\mu(U_{n,g})\geq 4.5$.
\end{lemma}

\begin{proof}
Define the numbers $a_g=2^k,\;a_i=a_{g-i}=a_{g+i}=2^{k-i}$ for $i\in\{1,2,\cdots,k\}$ and $a=\sqrt{\underset{j=1}{\overset{n}{\sum}}\;a_j^2}$.
We have
\begin{eqnarray}
a^2 & = & a_g^2 + 3(a_{g+1}^2+\cdots+a_{g+k-1}^2)+2a_n^2\nonumber\\
    & = & 2^{2k} + 3(2^{2k-2}+\cdots+2^2)+2\nonumber\\
    & = & 2(2^{2k}-1).\nonumber
\end{eqnarray}
Let $X=(x_1,\cdots,x_{n})$ be the unit vector, where $x_j=\frac{a_j}{a}$ for $j\in\{1,2,\cdots,n\}$. Then
\begin{eqnarray}
\underset{\{i,j\}\in E}{\sum}\;(x_i+x_j)^2 & = & 2\left(\underset{j=1}{\overset{n}{\sum}}\;x_j^2\right)+x_g^2-x_n^2+6\left(\underset{j=g}{\overset{n-1}\sum}\; x_jx_{j+1}\right)\nonumber\\
  & = & 2+ \frac{1}{a^2}(2^{2k}-1)+\frac{6}{a^2}(2^{2k-1}+2^{2k-3}+\cdots+2^3+2)\nonumber\\
  & = & 2+ \frac{1}{a^2}(2^{2k}-1)+\frac{4}{a^2}(2^{2k}-1)\nonumber\\
  & = & 4.5 \nonumber
\end{eqnarray}
So $\mu(U_{n,g})\geq X^{T}B(U_{n.g})X=\underset{\{i,j\}\in E}{\sum}\;(x_i+x_j)^2=4.5$. This completes the proof.
\end{proof}

\begin{corollary}\label{cor:4.5}
If $g$ is even and $n\geq \frac{3g}{2}$, then $\lambda(U_{n,g})\geq 4.5$.
\end{corollary}

\begin{proof}
This follows from Lemmas \ref{lem:bipartite} and \ref{lem:4.5}.
\end{proof}

We now consider the graph $U_{n,g}$ with $g$ odd. Let $\overline{U}_{n,g}$ be the graph obtained from $U_{n,g}$ by deleting the edge $\{\frac{g-1}{2},\frac{g+1}{2}\}$.
Then $\overline{U}_{n,g}$ is a tree, and so $\lambda(\overline{U}_{n,g})=\mu(\overline{U}_{n,g})$.

\begin{proposition}\label{prop:ung2}
Let $G=\overline{U}_{n,g}$ and $X=[x_1,\cdots,x_n]^T$ be a positive eigenvector of $B(G)$ corresponding to $\mu=\mu(G)$. Then the following hold:
\begin{enumerate}
\item[$(i)$] $x_j=x_{g-j}$ for $j\in\{1,\cdots, \frac{g-1}{2}\}$.
\item[$(ii)$] $x_j>x_{j+1}$ for $j\in\{g,\cdots, n-1$\}.
\item[$(iii)$] $x_g>x_1$ and $x_j>x_{j+1}$ for $j\in\{1,\cdots,\frac{g-3}{2}$\}.
\item[$(iv)$] If $\mu\geq 4.383$, then $x_i>x_{g+2i}$, where $1\leq i\leq\text{min }\{\frac{g-1}{2}, \lfloor\frac{n-g}{2}\rfloor\}$.
\end{enumerate}
\end{proposition}

\begin{proof}
The proof of $(i)$, $(ii)$ and $(iii)$ is similar to that of Proposition \ref{prop:ung1}$(i)$ and $(ii)$. We now prove $(iv)$. Along the path $g(g+1)(g+2)\cdots  n$ in $\overline{U}_{n,g}$, by Lemma \ref{lem:fx1}, we have
\begin{eqnarray}
x_g & = & f_{n-g}(\mu)x_{g+1}\nonumber\\
    & = & f_{n-g}(\mu)f_{n-g-1}(\mu)x_{g+2}\nonumber\\
    & \vdots & \nonumber\\
    & = & f_{n-g}(\mu)f_{n-g-1}(\mu)\cdots f_2(\mu)f_1(\mu)x_n.\nonumber
\end{eqnarray}
Similarly, along the path $g12\cdots\frac{g-1}{2}$ in $\overline{U}_{n,g}$, we have
$$x_g  = f_{\frac{g-1}{2}}(\mu)x_1  =  f_{\frac{g-1}{2}}(\mu)f_{\frac{g-3}{2}}(\mu)x_2=\cdots =f_{\frac{g-1}{2}}(\mu)\cdots f_1(\mu)x_{\frac{g-1}{2}}.$$
Thus, for $1\leq i\leq\text{min }\{\frac{g-1}{2}, \lfloor\frac{n-g}{2}\rfloor\}$, $$f_{\frac{g-1}{2}}(\mu)\cdots f_{\frac{g-2i+1}{2}}(\mu)x_{i}=x_g=f_{n-g}(\mu)\cdots f_{n-g-2i+1}(\mu)x_{g+2i}.$$
Since $\mu\geq 4.383$, $f_{n-g}(\mu)\cdots f_{n-g-2i+1}(\mu)>f_{\frac{g-1}{2}}(\mu)\cdots f_{\frac{g-2i+1}{2}}(\mu)$ by Lemma \ref{lem:fx}$(iii)$. So $x_i>x_{g+2i}$. This completes the proof.
\end{proof}

\begin{lemma}\label{rem1}
If $g\geq 5$ is odd and $n\geq g+2$, then $\mu(\overline{U}_{n,g})=\lambda(U_{n,g})>4.383$.
\end{lemma}

\begin{proof}
Since $g$ is odd, $\lambda(U_{n,g})=\lambda(\overline{U}_{n,g})$ by Proposition \ref{thm:guo}. Now consider the star graph $S$ on $4$ vertices. Form a new graph $S_1$ by appending one vertex to each of the pendant vertices of $S$. Using MATLAB, we have $\lambda(S_1)\approx 4.4142>4.383$. So for odd $g\geq 5$ and $n\geq g+2$,  $\mu(\overline{U}_{n,g})=\lambda(\overline{U}_{n,g})\geq \lambda(S_1)>4.383$  (see Lemma \ref{lem:bipartite}).
\end{proof}

\section{Proof of Theorem \ref{main}}

For a unicyclic graph $G$ on $n$ vertices with girth $g\geq 3$ ($n>g$), we take the vertices on the cycle of $G$ as $1,2,\cdots,g$ and denote by $l_i$ the number of vertices on the tree attached to the vertex $i$ for $1\leq i\leq g$. We define $C_G$ to be the set of all vertices $i$ on the cycle of $G$ for which $l_i\geq 1$. Then $l_1+\cdots+l_g=n-g$ and $1\leq |C_G|\leq g$.

For two distinct vertices $u$ and $v$ of $G$, we denote by $d(u,v)$ the {\it distance} between $u$ and $v$ (that is, the length of a shortest path between $u$ and $v$). The following lemma is useful for us.

\begin{lemma}\label{lem1:even}
Let $G$ be a unicyclic graph on $n$ vertices with girth $g$. Suppose that $|C_G|=r\geq 2$. Then the following hold:
\begin{enumerate}
\item[$(a)$] If $n\geq 2g-1$, then there exists two vertices $i,j\in C_G$ with $l_{i}\geq d(i,j)$.
\item[$(b)$] If $n\geq 3g-1$, then there exists two vertices $i, j\in C_G$ with $l_{i}\geq 2d(i,j)$.
\end{enumerate}
\end{lemma}

\begin{proof}
Let $i_1,i_2,\cdots,i_r$ be the vertices in $C_G$ with $i_1<i_2<\cdots<i_r$. Set $d_j=d(i_j,i_{j+1})$ for $1\leq j\leq r-1$ and $d_r=d(i_r,i_{1})$.

To prove $(a)$, it is enough to show that $d_j\leq \max\{l_j,l_{j+1}\}$ for some $1\leq j\leq r-1$ or $d_r\leq \max\{l_r,l_{1}\}$. Suppose that this is not true. Then $l_j\leq d_j-1$ for $1\leq j\leq r$. Now $g-1\leq l_1+\cdots+l_r\leq d_1+\cdots+d_r-r\leq g-r$. This implies that $r\leq 1$, a contradiction.

To prove $(b)$, we show that $2d_j\leq \max\{l_j,l_{j+1}\}$ for some $1\leq j\leq r-1$ or $2d_r\leq \max\{l_r,l_{1}\}$. If this is not true, then $l_j\leq 2d_j-1$ for $1\leq j\leq r$. Now $2g-1\leq l_1+\cdots+l_r\leq 2(d_1+\cdots+d_r)-r\leq 2g-r$. This gives $r\leq 1$, a contradiction.
\end{proof}

\subsection{The case $g$ even}

Let $G$ be a unicyclic graph on $n$ vertices with even girth $g$. By a repeated use of graph operations consisting of grafting of edges we can transform the graph $G$ into a new graph $G_1$ that has a path $P_i$ (on $l_i$ vertices) appended to the vertex $i$ for each $i\in C_G$ and that $\lambda(G_1)\leq \lambda(G)$. This is possible by Proposition \ref{thm:gkl}.

\begin{proposition}\label{prop1:even}
Let $G$ be a unicyclic graph on $n$ vertices with even girth $g$ which is not isomorphic to $U_{n,g}$. If $|C_G|=1$, then $\lambda(U_{n,g})< \lambda(G)$.
\end{proposition}

\begin{proof}
Since $G$ is not isomorphic to $U_{n,g}$, we use the operation grafting of edges at least once to get the graph $G_1$. Since $G$ is bipartite, $\lambda(G_1)< \lambda(G)$ by Proposition \ref{thm:gkl}. Now $|C_G|=1$ implies that $G_1$ is isomorphic to $U_{n,g}$. So $\lambda(U_{n,g})< \lambda(G)$.
\end{proof}

\begin{proposition}\label{prop2:even}
Let $G$ be a unicyclic graph on $n$ vertices with even girth $g$. Suppose that $|C_G|=r\geq 2$ and that we can arrange the vertices in $C_G$ in some ordering, say $i_1,i_2,\cdots,i_r$, such that $l_{i_1}+\cdots+l_{i_{j-1}}\geq d(i_{1},i_{j})$ for $2\leq j\leq r$. If $\lambda(U_{n,g})\geq 4.5$, then $\lambda(U_{n,g})< \lambda(G)$.
\end{proposition}

\begin{proof}
Consider the graph $G_1$ obtained from $G$ as above. For  $1\leq j\leq r$, let $i_j1$ be the vertex on the path $P_{i_j}$ adjacent to the vertex $i_j$, and  $i_jl_{i_j}$ be the pendant vertex of $P_{i_j}$. That is, the path $P_{i_j}$ is $(i_j1)(i_j2)\cdots(i_jl_{i_{j}})$, where the vertex $i_j1$ is adjacent to $i_j$. Note that $i_j1$ and $i_jl_{i_j}$ are the same vertices if $l_{i_j}=1$. With the sequence of vertices $i_1,i_2,\cdots,i_r$, we perform the following graph operations on $G_1$:
\begin{quote}
Delete the edge $\{i_2,i_21\}$ and add a new edge $\{i_1l_{i_1},i_21\}$, delete the edge $\{i_3,i_31\}$ and add a new edge $\{i_2l_{i_2},i_31\}$ and so on.
\end{quote}
Let $G_2$ be the new graph thus obtained from $G_1$. Then $G_2$ is isomorphic to $U_{n,g}$. Since $g$ is even, both $G_1$ and $G_2$ are bipartite. So, by Lemma \ref{lem:unitary}, $\mu(G_t)=\lambda(G_t)$ for $t=1,2$. Let $X$ be the positive unit eigenvector of $B(G_2)$ corresponding to $\mu(G_2)$. Then
$\lambda(G_1)-\lambda(G_2)=\mu(G_1)-\mu(G_2)\geq X^T B(G_1)X-\mu(G_2)=X^T B(G_1)X- X^T B(G_2)X=
[(x_{i_2}+x_{i_21})^2-(x_{i_1l_{i_1}}+x_{i_21})^2]+[(x_{i_3}+x_{i_31})^2-(x_{i_2l_{i_2}}+x_{i_31})^2]+\cdots+
[(x_{i_r}+x_{i_r1})^2-(x_{i_{r-1}l_{i_{r-1}}}+x_{i_r1})^2]$.
Since $\lambda(U_{n,g})\geq 4.5$ and $l_{i_1}+\cdots+l_{i_{j-1}}\geq d(i_{1},i_{j})$ for $2\leq j\leq r$, Proposition \ref{prop:ung1}$(v)$ implies that $x_{i_2}>x_{i_1l_{i_1}},x_{i_3}>x_{i_2l_{i_2}},\cdots,x_{i_r}>x_{i_{r-1}l_{i_{r-1}}}$. So $\lambda(G_1)-\lambda(G_2)>0$. Hence $\lambda(G)\geq \lambda(G_1)>\lambda(G_2)=\lambda(U_{n,g})$.
\end{proof}

As an immediate consequence of Propositions \ref{prop1:even} and \ref{prop2:even}, we have

\begin{corollary}\label{cor1:even}
Let $G$ be a unicyclic graph on $n$ vertices with even girth $g$ which is not isomorphic to $U_{n,g}$. If $l_k\geq \frac{g}{2}$ for some vertex $k$ in $C_G$, then $\lambda(U_{n,g})< \lambda(G)$.
\end{corollary}

\begin{proof}
Since $l_k\geq \frac{g}{2}$, $\lambda(U_{n,g})\geq 4.5$ by Corollary \ref{cor:4.5}. Now the result follows from Proposition \ref{prop1:even} for $|C_G|=1$ and from Proposition \ref{prop2:even} for $|C_G|\geq 2$. In the latter case, we can take any ordering of the vertices in $C_G$ starting with the vertex $k$. Note that the distance between two vertices in $C_G$ is at most $\frac{g}{2}$.
\end{proof}

\begin{proposition}\label{prop3:even}
Let $G$ be a unicyclic graph on $n$ vertices with even girth $g$. Suppose that $|C_G|=r\geq 2$ and that $n\geq 2g-1$. Then $\lambda(U_{n,g})< \lambda(G)$.
\end{proposition}

\begin{proof}
We have $\lambda(U_{n,g})\geq 4.5$ by Corollary \ref{cor:4.5}. So, by Proposition \ref{prop2:even}, it is enough to show that we can arrange the vertices in $C_G$ in some ordering $i_1,i_2,\cdots,i_r$ such that $l_{i_1}+\cdots+l_{i_{j-1}}\geq d(i_{1},i_{j})$ for $2\leq j\leq r$. We shall prove this by induction on $r$. By Lemma \ref{lem1:even}$(a)$, there exists two distinct vertices $i_1$ and $i_2$ in $C_G$ such that $l_{i_2}\geq d(i_1,i_2)$. If $r=2$, then we take the ordering $i_2,i_1$.

So assume that $r>2$. Now consider the graph $G_1$ obtained from $G$ as above. Let $K$ be the graph obtained from $G_1$ by disconnecting the path $P_{i_1}$ from the vertex $i_1$ and appending it to the pendant vertex of the path $P_{i_2}$. We have $C_G=C_{G_1}$ and $C_K=C_G-\{i_1\}$. Set $l_{i_2}'=l_{i_1}+l_{i_2}$. Since $|C_K|<r$, applying induction hypothesis to the graph $K$, we can get an ordering $k_2,\cdots,k_t(=i_2),\cdots,k_r$ of the vertices in $C_K$ such that $l_{k_2}+\cdots+l_{k_{j-1}}\geq d(k_{2},k_{j})$ for $3\leq j\leq t$, $l_{k_2}+\cdots+l_{k_{t-1}}+l_{k_t}'\geq d(k_{2},k_{t+1})$ and  $l_{k_2}+\cdots+l_{k_{t-1}}+l_{k_t}'+l_{k_{t+1}}+\cdots+l_{k_{j-1}}\geq d(k_{2},k_{j})$ for $t+2\leq j\leq r$. Now the facts $l_{i_2}\geq d(i_1,i_2)$ and $l_{i_2}'=l_{i_1}+l_{i_2}$ imply that the ordering $k_2,\cdots,k_t(=i_2),i_1,k_{t+1},\cdots,k_r$ of the vertices in $C_G$ satisfies our requirement.
\end{proof}

Now Theorem \ref{main}(1) follows from Propositions \ref{prop1:even} and \ref{prop3:even}.\\

Consider the case $g=4$. There is only one unicyclic graph on five vertices. For $n\geq 6$, we have $\lambda(U_{n,4})\geq 4.5$ by Corollary \ref{cor:4.5}. By Theorem \ref{main}(1), $U_{n,4}$ uniquely minimizes the Laplacian spectral radius when $n\geq 7$. There are four non-isomorphic unicyclic graphs on six vertices and exactly one of them does not satisfy the conditions in Propositions \ref{prop1:even} or \ref{prop2:even}. This graph $G$ is obtained from a cycle of length four by appending one pendant vertex to each of its two opposite vertices. Using MATLAB we can see that $\lambda(G)\approx 4.73205 > 4.5615\approx\lambda(U_{6,4})$.

Now, consider the case $g=6$. There is only one unicyclic graph when $n=7$ and five non-isomorphic unicyclic graphs when $n=8$. In the latter case, we can verify that $U_{8,6}$ (with $\lambda(U_{8,6})\approx 4.4989$) uniquely minimizes the Laplacian spectral radius. For $n\geq 9$, we have $\lambda(U_{n,6})\geq 4.5$ by Corollary \ref{cor:4.5}. By Theorem \ref{main}(1), $U_{n,6}$ uniquely minimizes the Laplacian spectral radius for $n\geq 11$. There are three graphs for $n=10$ (and also for $n=9$) with $|C_G|\geq 2$ for which the hypothesis in Proposition \ref{prop2:even} is not satisfied. In both cases, we can check that $U_{n,6}$ uniquely minimizes the Laplacian spectral radius. Thus we have the following.

\begin{proposition}\label{prop4:girth4-6}
For $g\in\{4,6\}$, the graph $U_{n,g}$ uniquely minimizes the Laplacian spectral radius over all unicyclic graphs on $n$ vertices with girth $g$.
\end{proposition}

Having seen that $U_{n,g}$ uniquely minimizes the Laplacian spectral radius when $g=4\text{ and }6$, it is natural to expect that the same might be true for all even $g$. However, this is false when $n$ is not large comparing to $g$. We give an example below.

\begin{example}
Consider a cycle on $10$ vertices. Let $v$ and $w$ be two opposite vertices on this circle, that is, $d(v,w)=5$. Add one pendant vertex to each of the vertices $v$ and $w$. The new graph is a unicyclic graph on $12$ vertices with girth $10$, denote it by $C_{10}^{1,1}$. We have $4.4383\approx\lambda(C_{10}^{1,1}) < \lambda(U_{12,10})\approx 4.4763$.
\end{example}

\subsection{The case $g$ odd}

\begin{proposition}\label{prop1:odd}
Let $H$ be a unicyclic graph on $n$ vertices with odd girth $g$ which is not isomorphic to $U_{n,g}$. If $|C_H|=1$, then $\lambda(U_{n,g})< \lambda(H)$.
\end{proposition}

\begin{proof}
Let $C_H=\{j\}$ and $H_1$ be the bipartite graph obtained from $H$ by deleting the edge opposite to the vertex $j$. Since $H$ is not isomorphic to $U_{n,g}$, we can use the operation grafting of edges on the graph $H_1$ (at least once) to get a new graph $H_2$ which is isomorphic to $\overline{U}_{n,g}$. By Proposition \ref{thm:gkl}, $\lambda(H_2)< \lambda(H_1)$. Now Proposition \ref{thm:guo} and Lemma \ref{lem:courant} imply that $\lambda(U_{n,g})=\lambda(\overline{U}_{n,g})=\lambda(H_2)< \lambda(H_1)\leq \lambda(H)$.
\end{proof}

Consider the case $g=3$. Let $H$ be a unicyclic graph on $n$ vertices with girth three which is not isomorphic to $U_{n,3}$. If $|C_H|=1$, then $\lambda(U_{n,3})< \lambda(H)$ by Proposition \ref{prop1:odd}. Assume that $|C_H|\geq 2$.
For $n=5$, we can see using MATLAB that $U_{5,3}$ uniquely minimizes the Laplacian spectral radius. So assume that $n\geq 6$. Choose an edge $\{u, v\}$ on
the cycle of $H$ such that $H_1 = H-\{u,v\}$ is not a path. Then $H_1$ is a tree. By a finite sequence of grafting of edges (at least once) we can transform $H_1$ to a new tree $H_2$ such that $H_2$ is isomorphic to $\overline{U}_{n,3}$ and that $\lambda(H_2)<\lambda(H_1)$ (see Proposition \ref{thm:gkl}). Now, by Proposition \ref{thm:guo} and Lemma \ref{lem:courant}, $\lambda(U_{n,3})=\lambda(\overline{U}_{n,3})=\lambda(H_2)<\lambda(H_1)\leq\lambda(H)$. Thus we have the following.

\begin{proposition}\label{prop:girth3}
The graph $U_{n,3}$ uniquely minimizes the Laplacian spectral radius over all unicyclic graphs on $n$ vertices with girth three.
\end{proposition}

We give an example where $U_{n,g}$ does not minimize the Laplacian spectral radius among all unicyclic graphs on $n$ vertices with odd girth $g$.

\begin{example}
Consider a cycle on $7$ vertices. Let $v$ and $w$ be two vertices on this circle with $d(v,w)=3$. Add one pendant vertex to each of the vertices $v$ and $w$. The new graph is a unicyclic graph on $9$ vertices with girth $7$, denote it by $C_{7}^{1,1}$. We have $\lambda(C_{7}^{1,1})\approx 4.4142 < 4.4605\approx\lambda(U_{9,7})$.
\end{example}

Now, let $H$ be a unicyclic graph on $n$ vertices with odd girth $g\geq 5$. By a sequence of grafting of edges we can transform $H$ into a new graph $H_1$ that has a path $P_i$ (on $l_i$ vertices) appended to the vertex $i$ for each $i\in C_H$ and that $\lambda(H_1)\leq \lambda(H)$ (see Proposition \ref{thm:gkl}). 

The proof of the next result is similar to that of Proposition \ref{prop2:even} with some modifications. We write the proof with necessary changes for the sake of completeness.

\begin{proposition}\label{prop2:odd}
Let $H$ be a unicyclic graph on $n$ vertices with odd girth $g\geq 5$. Suppose that $|C_H|=r\geq 2$ and that we can arrange the vertices in $C_H$ in some ordering, say $i_1,i_2,\cdots,i_r$, such that $l_{i_1}+\cdots+l_{i_{j-1}}\geq 2d(i_{1},i_{j})$ for $2\leq j\leq r$. Then $\lambda(U_{n,g})< \lambda(H)$.
\end{proposition}

\begin{proof}
Consider the graph $H_1$ obtained from $H$ as above. For  $1\leq j\leq r$, let $i_j1$ be the vertex on the path $P_{i_j}$ adjacent to the vertex $i_j$, and  $i_jl_{i_j}$ be the pendant vertex of $P_{i_j}$. Let $H_2$ be the graph obtained from $H_1$ by deleting the edge opposite to the vertex $i_1$. With the given ordering of vertices $i_1,i_2,\cdots,i_r$, perform the same sequence of graph operations on $H_2$ as we have done on the graph $G_1$ in the proof of Proposition \ref{prop2:even} to get a new graph $H_3$. Then $H_3$ is isomorphic to $\overline{U}_{n,g}$. Since $H_2$ and $H_3$ are bipartite, $\mu(H_t)=\lambda(H_t)$ for $t=2,3$. Let $X$ be the positive unit eigenvector of $B(H_3)$ corresponding to $\mu(H_3)$. Then
$\lambda(H_2)-\lambda(H_3)=\mu(H_2)-\mu(H_3)\geq X^T B(H_2)X-\mu(H_3)=X^T B(H_2)X- X^T B(H_3)X=
[(x_{i_2}+x_{i_21})^2-(x_{i_1l_{i_1}}+x_{i_21})^2]+[(x_{i_3}+x_{i_31})^2-(x_{i_2l_{i_2}}+x_{i_31})^2]+\cdots+
[(x_{i_r}+x_{i_r1})^2-(x_{i_{r-1}l_{i_{r-1}}}+x_{i_r1})^2]$.
Since $\mu(\overline{U}_{n,g})> 4.383$ (Lemma \ref{rem1}) and $l_{i_1}+\cdots+l_{i_{j-1}}\geq 2d(i_{1},i_{j})$ for $2\leq j\leq r$, Proposition \ref{prop:ung2}$(iv)$ implies that $\lambda(H_2)-\lambda(H_3)>0$. So $\lambda(U_{n,g})=\overline{U}_{n,g})=\lambda(H_3)<\lambda(H_2)\leq \lambda(H)$.
\end{proof}

As a consequence of Propositions \ref{prop1:odd} and \ref{prop2:odd}, we have

\begin{corollary}\label{cor1:odd}
Let $H$ be a unicyclic graph on $n$ vertices with odd girth $g\geq 5$ which is not isomorphic to $U_{n,g}$. If $l_k\geq g-1$ for some vertex $k$ in $C_H$, then $\lambda(U_{n,g})< \lambda(H)$.
\end{corollary}

\begin{proof}
This follows from Proposition \ref{prop1:odd} for $|C_H|=1$ and from Proposition \ref{prop2:odd} for $|C_H|\geq 2$. In the latter case, we can take any ordering of the vertices in $C_H$ starting with the vertex $k$.
\end{proof}

\begin{proposition}\label{prop3:odd}
Let $H$ be a unicyclic graph on $n$ vertices with odd girth $g\geq 5$. Suppose that $|C_H|=r\geq 2$ and that $n\geq 3g-1$. Then $\lambda(U_{n,g})< \lambda(H)$.
\end{proposition}

\begin{proof}
The proof of this result is similar to that of Proposition \ref{prop3:even} (In the proof, one has to replace Proposition \ref{prop2:even} by Proposition \ref{prop2:odd}, Lemma \ref{lem1:even}$(a)$ by Lemma \ref{lem1:even}$(b)$, $G$ by $H$ and $G_1$ by $H_1$).
\end{proof}

Now Theorem \ref{main}(2) follows from Propositions \ref{prop1:odd}, \ref{prop:girth3} and \ref{prop3:odd}.

\end{document}